\documentclass[doublespacing]{elsart}
\usepackage{natbib,amsmath,amssymb,hyperref,enumerate}
\def\alephe{S}

\begin{document}

\begin{frontmatter}

\title{On the existence of some ARCH($\infty$) processes}

\author[randal]{Randal Douc} \author[francois]{Fran\c{c}ois
  Roueff} \author[philippe]{Philippe Soulier\corauthref{cor1}}
\ead{philippe.soulier@u-paris10.fr} \corauth[cor1]{Corresponding
  author} 
 \address[randal]{Ecole
  Polytechnique, CMAP, 91128 Palaiseau Cedex, France}
\address[francois]{Telecom Paris, CNRS LTCI, 75634 Paris Cedex 13,
  france}
\address[philippe]{Universit\'e Paris X, Laboratoire
  MODAL'X, 92000 Nanterre, France}

\begin{abstract}
  A new sufficient condition for the existence of a stationary causal
  solution of an ARCH($\infty$) equation is provided. This condition
  allows to consider coefficients with power-law decay, so that it
  can be applied to the so-called FIGARCH processes, whose existence
  is thus proved.
\end{abstract}


\end{frontmatter}

\section{Introduction}
It can arguably be said that autoregressive conditionnally
heteroskedastic (ARCH) and long memory processes are two success
stories of the nineties, so that they were bound to meet. Their
tentative offspring was the FIGARCH process, introduced by
\cite{baillie:bollerslev:mikkelsen:1996} without proving its
existence, which has remained controversial up to now. More precisely,
the FIGARCH($p,d,q)$ process is the solution of the equations
\begin{align}
  X_n & = \sigma_n z_n \; ,  \label{eq:levels} \\
  \sigma_n^2 & = a_0 + \left\{ I-(I-L)^d \; \frac{\theta(L)}{\phi(L)} \right\} X_n^2 \;
  , \label{eq:figarch0d0}
\end{align}
where $\{z_n\}$ is an i.i.d. sequence with zero mean and unit
variance, $a_0>0$, $d \in(0,1)$, $L$ is the backshift operator and $(I-L)^d$
is the fractional differencing operator:
\begin{align*}
  (I-L)^d = I + \sum_{j=1}^\infty  \frac{(-d)(1-d)\cdots(j-1-d)}{j!}
  L^j \; ,
\end{align*}
and $\theta$ and $\phi$ are polynomials such that $\theta(0)=\phi(0) = 1$,
$\phi(z)\neq0$ for all complex number $z$ in the closed unit disk and
the coefficients of the series expansion of $1-(1-z)^d\theta(z)/\phi(z)$ are
nonnegative. Then the coefficients $\{a_j\}_{j\geq1}$ defined by
$\sum_{j=1}^\infty a_jL^j=I-(I-L)^d\theta(L)/\phi(L)$ satisfy $a_j \sim c
j^{-d-1}$ for some constant $c>0$ and $\sum_{j=1}^\infty a_j=1$.

These processes are subcases of what can be called IARCH($\infty)$,
defined as solutions of the equations~(\ref{eq:levels}) and
\begin{align} \label{eq:archinfty}
  \sigma_n^2 & = a_0 + \sum_{j=1}^\infty a_j X_{n-j}^2 \; ,
\end{align}
for some sequence $\{a_j\}$ such that $a_0>0$ and $\sum_{j=1}^\infty
a_j = 1$.  The letter $I$ stands for integrated, by analogy to ARIMA
processes. An important property of such processes is that a
stationary solution necessarily has infinite variance.  Indeed, if
$\sigma^2=\mathbb{E}[\sigma_n^2]<\infty$, then $\mathbb{E}[X_n^2] =
\sigma^2$ and~(\ref{eq:archinfty}) implies $\sigma^2 = a_0 +
\sigma^2$, wich is impossible.  If the condition $\sum_{j=1}^\infty
a_j = 1$ is not imposed, a solution to equations~(\ref{eq:levels})
and~(\ref{eq:archinfty}) is simply called an ARCH($\infty$) process.

A solution of an ARCH($\infty)$ equation is said to be causal with
respect to the i.i.d. sequence $\{z_n\}$ if for all $n$, $\sigma_n$ is
$\mathcal{F}_{n-1}^z$ measurable, where $\mathcal{F}_n^z$ is the
sigma-field generated by $\{z_n, z_{n-1}, \ldots\}$. Note that to
avoid trivialities, here and in the following, $\sigma_n$ is the positive
square root of $\sigma_n^2$.  There exists an important literature on
ARCH($\infty$), IARCH($\infty$) and FIGARCH processes. For a recent
review, see for instance
\cite{giraitis:leipus:surgailis:2007}.  The known conditions for the
existence of stationary causal conditions to ARCH equations are always
a compromise between conditions on the distribution of the innovation
sequence $\{z_n\}$ and summability conditions on the coefficients
$\{a_j, j\geq1\}$. \cite{giraitis:surgailis:2002} provides a necessary
and sufficient condition for the solution to have finite fourth
moment. The only rigorous result in the IARCH($\infty$) case was
obtained by \cite{kaza:leipus:2003}. They prove the existence of a
causal stationary solution under the condition that the coefficients
$a_j$ decay geometrically fast, which rules out FIGARCH processes, and
on a mild condition on the distribution of $z_0$.

The purpose of this paper is to provide a new sufficient condition for
the existence of a stationary solution to an ARCH($\infty$) equation,
which allows power-law decay of the coefficients $a_j$s, even in the
IARCH($\infty$) case. This condition is stated in
Section~\ref{sec:sufficient}. It is applied to the IARCH($\infty$)
case in Section~\ref{sec:iarch} and the existence of a stationary
solution to the FIGARCH equation is proved. Further research
directions are given in Section~\ref{sec:problems}. In particular, the
memory properties of FIGARCH processes are still to be investigated.
This is an important issue, since the original motivation of these
processes was the modelling of long memory in volatility.

\section{A sufficient condition for the existence of ARCH($\infty$) processes} 
\label{sec:sufficient}

\begin{thm} \label{thm:csexistence}
  Let $\{a_j\}_{j\geq0}$ be a sequence of nonnegative real numbers and
  $\{z_k\}_{k\in\Zset}$ a sequence of i.i.d. random variables. For $p>0$, define
  $$
  A_p = \sum_{j=1}^\infty a_j^p
  \quad\text{and}\quad
  \mu_p = \mathbb{E}[z_0^{2p}]\;. 
  $$
If there exists $p\in(0,1]$ such that  
  \begin{gather} \label{eq:cs}
    A_p \mu_p < 1, 
  \end{gather}
  then there exists a  strictly stationary solution of the ARCH($\infty$)
  equation:
  \begin{align}
\label{eq:EqARCH}
X_n & = \sigma_n z_n \; , \\
\sigma_n^2 & = a_0 + \sum_{j=1}^\infty a_j X_{n-j}^2 \; , \label{eq:volatilite}
  \end{align} 
given by~(\ref{eq:EqARCH}) and 
\begin{gather}
  \sigma_n^2 = a_0 + a_0 \sum_{k=1}^\infty \sum_{j_1,\dots,j_k\geq1}
  a_{j_1}\dots a_{j_k} z_{n-j_1}^2 \dots z_{n-j_1-\dots-j_k}^2 \; .
  \label{eq:solarch}
\end{gather}
The process $\{X_n\}$ so defined is the unique
causal stationary solution
to equations~(\ref{eq:EqARCH}) and~(\ref{eq:volatilite}) such that
$\mathbb{E}[|X_n|^{2p}]<\infty$.
\end{thm}

\noindent{\em Proof.}
Denote $\xi_k=z_k^2$, so that $\mathbb{E}[\xi_k^p]=\mu_p$, and define the
$[0,\infty]$-valued r.v.
  \begin{gather} \label{eq:aleph}
  \alephe_0 = a_0 + a_0 \sum_{k=1}^\infty \sum_{j_1,\dots,j_k\geq1} a_{j_1}\dots a_{j_k}
  \xi_{-j_1} \dots \xi_{-j_1-\dots-j_k}
\end{gather}
Since $p \in (0,1]$, we apply the inequality $(a+b)^p \leq a^p + b^p$
valid for all $a,b\geq0$ to $S_0^p$:
\begin{align*}
  \alephe_0^p & \leq a_0^p + a_0^p \sum_{k=1}^\infty
  \sum_{j_1,\dots,j_k\geq1} a_{j_1}^p\dots a_{j_k}^p
\xi_{-j_1}^p \dots \xi_{-j_1-\dots-j_k}^p \; .
\end{align*}
Then, by independence of the $\xi_j$'s, we obtain
\begin{align}\nonumber
  \mathbb{E}[\alephe_0^p] & \leq a_0^p + a_0^p \sum_{k=1}^\infty
  \sum_{j_1,\dots,j_k\geq1} a_{j_1}^p\dots a_{j_k}^p
  \mathbb{E}[\xi_{-j_1}^p \dots
  \xi_{-j_1-\dots-j_k}^p] \\
  \label{eq:alephBound}
  & = a_0^p\left[1+\sum_{k=1}^\infty (\mu_pA_p)^k\right] =
  \frac{a_0^p}{1 - A_p \mu_p} \; ,
\end{align}
where we used~(\ref{eq:cs}). This bound shows that $\alephe_0<\infty$
a.s. and the sequence
$$
\alephe_n = a_0 + a_0 \sum_{k=1}^\infty \sum_{j_1,\dots,j_k\geq1}
a_{j_1}\dots a_{j_k} \xi_{n-j_1} \dots \xi_{n-j_1-\dots-j_k},\quad
n\in\Zset\;,
$$
is a sequence of a.s. finite r.v.'s. Since only nonnegative numbers
are involved in the summation, we may write
\begin{multline*}
  \sum_{j=1}^\infty a_j \alephe_{n-j}\xi_{n-j}
  =a_0\sum_{j_0=1}^\infty a_{j_0} \xi_{n-j_0}\\
  +a_0\sum_{j_0=1}^\infty a_{j_0}\xi_{n-j_0} \sum_{k=1}^\infty
  \sum_{j_1,\dots,j_k\geq1} a_{j_1}\dots a_{j_k}
  \xi_{n-j_0-j_1} \dots \xi_{n-j_0-j_1-\dots-j_k} \\
  =a_0\sum_{k=0}^\infty\sum_{j_0,j_1,\dots,j_k\geq1}a_{j_0}\dots
  a_{j_k}\xi_{n-j_0} \dots \xi_{n-j_0-j_1-\dots-j_k} \;.
\end{multline*}
Hence $\{\alephe_n,\,n\in\Zset\}$ satisfies the recurrence equation
$$
\alephe_n= a_0+\sum_{j=1}^\infty a_j \alephe_{n-j}\xi_{n-j} \; .
$$
The technique of infinite chaotic expansions used here is standard; it was
already used in the proof of \cite[Theorem 2.1]{kokoszka:leipus:2000}.  This
proves the existence of a strictly stationary solution for~(\ref{eq:EqARCH})
and~(\ref{eq:volatilite}) by setting $\sigma_n^2=\alephe_n$ and
$X_n=\sigma_nz_n$.  Using~(\ref{eq:alephBound}), we moreover have
$\mathbb{E}[|X_n|^{2p}]\leq \mu_pa_0^p/(1 - A_p \mu_p)$.

Suppose now that $\{X_n\}$ is a strictly stationary causal solutions of the
ARCH($\infty$) equations~(\ref{eq:EqARCH}) and~(\ref{eq:volatilite}). Then, for
any $q \geq 1$, the following expansion holds:
\begin{align} \label{eq:series}
    \sigma_n^2 & = a_0 + a_0 \sum_{k=0}^q \sum_{j_1,\dots,j_k\geq1}a_{j_1}\dots
  a_{j_k} \xi_{n-j_1} \dots \xi_{n-j_1-\dots-j_k} \\
  & + \sum_{j_1,\dots,j_{q+1}\geq1} a_{j_1} \dots a_{j_{q+1}} \xi_{n-j_1} \dots
  \xi_{n-j_1-\dots-j_q} X_{n-j_1-\cdots-j_{q+1}}^2 \; .\label{eq:reste}
\end{align}
The last display implies that the series on the right-hand side of~(\ref{eq:series})
converges to $\alephe_n$ as $q\to\infty$. Denote by $R_{n,q}$ the
remainder term in~(\ref{eq:reste}). Since $\{X_n\}$ is a causal solution,
$X_{n-j_1-\cdots-j_{q+1}}$ is independent of $ \xi_{n-j_1} \dots
\xi_{n-j_1-\dots-j_q}$ for all $j_1,\dots,j_{q+1} \geq1$. Hence, for any $p\leq
1$,
\begin{align*}
  \mathbb{E}[R_{n,q}^p] \leq (A_p \mu_p)^q \mathbb{E}[X_0^{2p}] \; .
\end{align*}
If Assumption~(\ref{eq:cs}) holds and $\mathbb{E}[X_0^{2p}] < \infty$, then
$\mathbb{E}[\sum_{q\geq1}R_{n,q}^p]<\infty$ so that, as $q\to\infty$, $R_{n,q}\to0$ a.s., implying  $\sigma_n^2 = \alephe_n$
a.s.  \hfill$\Box$

\section{IARCH($\infty$)  processes}
\label{sec:iarch}
IARCH (Integrated ARCH) processes are particular ARCH($\infty$)
processes for which $A_1\mu_1=1$, or, equivalently up to a scale
factor,
\begin{gather} \label{eq:figarchCondition}
A_1=1 \quad\text{and}\quad \mu_1=1
\end{gather}
To the best of our knowledge, the only rigorous general result on
IARCH($\infty$) processes was obtained by \cite{kaza:leipus:2003}.
See \cite{giraitis:leipus:surgailis:2007} for a recent review.  In Theorem~2.1 of 
\cite{kaza:leipus:2003}, it is proved that if
\begin{gather} \label{eq:kazakC1}
  \mathbb{E}[|\log (z_0)|^2]<\infty \; ,\\
  \label{eq:kazakC2}
  \sum_i a_i q^i <\infty\quad\text{for some}\quad q>1 \; ,
\end{gather}
hold, then there exists a unique stationary causal solution to the
ARCH($\infty$) equations~(\ref{eq:EqARCH})-(\ref{eq:volatilite}).
Condition~(\ref{eq:kazakC1}) on the distribution of $z_0$ is mild, but
the condition~(\ref{eq:kazakC2}) rules out power-law decay of the
coefficients $\{a_j\}$. 

Theorem~\ref{thm:csexistence} yields the following sufficient condition
for the existence of a IARCH($\infty$) process.

\begin{cor} \label{cor:CNScond} If $A_1=1$ and $\mu_1=1$, (\ref{eq:cs}) holds
  for some $p\in(0,1]$ if and only if there exists $p^*<1$ such that $A_{p_*} <
  \infty$ and
  \begin{gather} 
    \sum_{j=1}^\infty a_i \log(a_i) + \mathbb{E}[z_0^2\log(z_0^2)] \in(0,\infty] \; .
    \label{eq:condition}
  \end{gather}
  Then, the process defined by~(\ref{eq:EqARCH}) and~(\ref{eq:solarch}) is a
  solution of the ARCH($\infty)$ equation and  $\mathbb{E}
  [|X_n|^{q}]<\infty$ for all $q\in[0,2)$ and $\mathbb{E}[X_n^2]=\infty$.
\end{cor}

\noindent{\em Proof.}
Since $a_i \leq 1$ for all $i\geq1$, it holds that $\sum_{j=1}^\infty a_i
\log(a_i)\leq0$ and the convexity of the function $x\mapsto x\log(x)$ implies
$\mathbb{E}[z_0^2\log(z_0^2)]\geq0$.

First assume that there exists $p\in(0,1]$ such that~(\ref{eq:cs}) holds. Since
$A_1=\mu_1=1$, then necessarily, $p<1$ and for all $q\in[p,1]$, $A_q<\infty$.
Thus we can define the function $\phi:[p,1]\to\Rset$ by
$$
\phi(q)=\log(A_q\mu_q)=\log\sum_{j=1}^\infty a_j^q+\log \mathbb{E}[z_0^{2q}]
\;.
$$
H{\"o}lder inequality implies that the functions
$q\mapsto\log\sum_{j=1}^\infty a_j^q$ and $q\mapsto\log
\mathbb{E}[z_0^{2q}]$ are both convex on $[p,1]$. Thus $\phi$ is also convex
on $[p,1]$ and, since $\phi(p)<0$ and $\phi(1)=0$, the left derivative
of $\phi$ at $1$, which is given by the left-hand side
of~(\ref{eq:condition}), is positive (possibly infinite).

Conversely suppose that there exists $p^\ast<1$ such that
$A_{p^\ast}<\infty$ and that~(\ref{eq:condition}) holds.  Then $\phi$
is a convex function on $[p^\ast,1]$ and~(\ref{eq:condition}) implies
that $\phi(q)<0$ for $q<1$ sufficiently close to 1.

By convexity of $\phi$ and since $\phi(1)=0$, we also get that
$A_p\mu_p<1$ implies $A_q\mu_q<1$ for all $q\in[p,1)$. Then, by
Theorem~\ref{thm:csexistence}, the process $\{X_n,\,n\in\Zset\}$
defined by~(\ref{eq:solarch}) and~(\ref{eq:volatilite}) is a solution
to the ARCH($\infty$) equation and satisfies $\mathbb{E}[|X_0|^q]<\infty$
for all positive $q<2$. \hfill $\Box$.

\noindent{\em Comments on Corollary~\ref{cor:CNScond}.}
\begin{enumerate}[(i)]
\item Condition~(\ref{eq:condition}) is not easily comparable to
  conditions~~(\ref{eq:kazakC1}) and~(\ref{eq:kazakC2}) of
  \cite{kaza:leipus:2003}.
  Condition~(\ref{eq:condition}) is not necessary to prove the existence of a
  causal stationary solution if the coefficients $a_j$ decay geometrically fast
  (in particular if there are only finitely many nonvanishing coefficients), as
  a consequence of \cite[Theorem 2.1]{kaza:leipus:2003}; however, this result
  does not prove that any moments of $X_n$ are finite, contrary to
  Corollary~\ref{cor:CNScond}.
\item It might also be of interest to note that the Lyapounov exponent
  of the FIGARCH process as defined in \cite{kaza:leipus:2003} is
  zero. So our result proves that such a feature is not in
  contradiction with strict stationarity.

\item In the specific case of IGARCH processes, which are particular parametric
  subclasses of IARCH($\infty$) processes, \cite{bougerol:picard:1992} have a
  different set of assumptions on the distribution of $z_0$: they assume that
  $\mathbb{P}(z_0^2=0)=0$ and that the support of the distribution of $z_0^2$
  is unbounded.
\item The moment $\mathbb{E}[z_0^2 \log(z_0^2)]$ can be arbitrarily
  large (possibly infinite) if the distribution of $z_0^2$ has a
  sufficiently heavy tail. It is infinite for instance if the
  distribution of $z_0^2$ is absolutely continuous with a density
  bounded from below by $1/(x^2\log^2(x))$ for $x$ large enough.
  In that case, condition~(\ref{eq:condition}) holds for any sequence $\{a_j\}$
  such that $A_{p^*}<\infty$ for some $p^*<1$.  This conditions allows for a
  power-law decay of the coeffficients $a_j$, for instance $a_j \sim
  cj^{-\delta}$, for some $\delta>1$.

\end{enumerate}

Corollary~\ref{cor:CNScond} can be used to prove the existence of a
causal strictly stationary solution to some FIGARCH($p,d,q$)
equations.  Let us illustrate this in the case of the FIGARCH($0,d,0$)
equation, that is~(\ref{eq:EqARCH}) and~(\ref{eq:volatilite}) with
$d\in(0,1)$, $a_0 > 0$ and $a_j=\pi_j(d)$ for all $j\geq1$, where
\begin{gather*}
  \pi_1(d) = d\; , \ \pi_j(d) = \frac{d(1-d)\cdots(j-1-d)}{j!} \; ,
  j\geq2 \; .
\end{gather*}

\begin{cor}
  Assume that $\{z_k\}_{k\in\Zset}$ a sequence of i.i.d. random
  variables, such that $\mathbb{E}[z_0^2]=1$ and
  $\mathbb{P}\{|z_0|=1\}<1$. Then there exists $d^*\in[0,1)$ such
  that, for all $d\in(d^*,1)$, the FIGARCH($0,d,0$) equation has a
  unique causal stationary solution satisfying
  $\mathbb{E}[|X_n|^{2p}]<\infty$ for all $p<1$.
\end{cor}
\noindent{\em Proof. }
For $d \in(0,1]$ and $p\in(1/(d+1),1]$, denote
\begin{gather*}
  H(p,d) = \log \sum_{j=1}^\infty \pi_j^p(d) \; , \ \ L(d) =
  \sum_{j=1}^\infty \pi_j(d) \log(\pi_j(d) ) \; .
\end{gather*}
For $d\in(0,1)$, $\pi_j(d) \sim c j^{-d-1}$, so that $H(p,d)$ is
defined on $(1/(d+1),1]$. Moreover, it is decreasing and convex with
respect to $p$, $H(1,d)=0$ and $\partial_p H(1,d) = L(d)$.  Also,
$\pi_j(d)/d$ is a decreasing function of $d$ and $\lim_{d\to1}\pi_j(d)
= 0$ for all $j\geq2$. Thus, by bounded (and monotone) convergence,
for all $p\in(1/2,1)$, it holds that $\lim_{d\to1, d<1} H(p,d)= 0$.
By convexity of $H$ with respect to $p$, the following bound holds:
\begin{gather*}
  0 \leq - L(d) \leq \frac{H(p,d)}{1-p} \; .
\end{gather*}
Hence $\lim_{d\to1} L(d) = 0$.  
By assumption, we have $\mathbb{E}[z_0^2\log(z_0^2)]>0$. 
This implies that there exists $d^*\in(0,1)$ such that $L(d)+
\mathbb{E}[z_0^2 \log(z_0^2)] >0$ (i.e.~(\ref{eq:condition}) holds) if
$d>d^*$. Thus Corollary~\ref{cor:CNScond} proves the existence of the
corresponding FIGARCH($0,d,0$) processes. \hfill $\Box$.

\noindent{\em Remark.}
It easily seen that $L(d) \leq \log(d)$ so that $\lim_{d\to0} L(d) =
-\infty$, i.e.~(\ref{eq:condition}) does not hold for small $d$. We
conjecture, but could not prove, that $L(d)$ is increasing, so
that~(\ref{eq:condition}) holds if \emph{and only if} $d>d^*$ (with
$d^*=0$ if $\mathbb{E}[\xi_0\log(\xi_0)]=\infty$).  But this does not
prove that the FIGARCH($0,d,0$) does not exist for $d\leq d^*$.

\section{Open problems}
\label{sec:problems}
Now that a proof of existence of some FIGARCH and related processes is
obtained under certain conditions, there still remain some open
questions. We state a few of them here.
\begin{enumerate}[(i)]
\item Condition~(\ref{eq:condition}) is not necessary for the
  existence of a stationary causal solution, but it implies finiteness
  of all moments up to 1 of $X_n^2$ (with of course $\mathbb{E}[X_n^2]
  = \infty$). The problem remains open to know if there exist a
  stationary solution under a mild assumption on $z_0$, such
  as~(\ref{eq:kazakC1}) for instance.  If a solution exists, say
  $\{X_n\}$, then, as seen in the proof of
  Theorem~\ref{thm:csexistence}, the sequence $\{\alephe_n\}$ defined
  in~(\ref{eq:aleph}) is well defined and $Y_n = \alephe_n^{1/2} z_n$
  is also a stationary causal solution wich satisfies moreover $Y_n^2
  \leq X_n^2$. But we cannot prove without more assumptions that these
  solutions are equal.
\item Tail behaviour of the marginal distribution of GARCH processes
  have been investigated by \cite{basrak:davis:mikosch:2002},
  following \cite{nelson:1990}, but there are no such results in the
  ARCH($\infty$) case.  Under suitable conditions, we have shown that
  the squares of the FIGARCH process $X_n^2$ have finite moments of
  all order $p<1$, but necessarily, $\mathbb{E}[X_n^2]=\infty$. Thus,
  it is natural to conjecture that perhaps under additional conditions
  on the distribution of $z_0$, the function $x\to\mathbb{P}(X_n^2>x)$
  is regularly varying with index -1.
\item The memory properties of the FIGARCH process are of course of
  great interest. The sequence $\{X_n\}$ is a strictly stationary
  martingale increment sequence, but $\mathbb{E}[X_n^2]=\infty$. So does it
  hold that the partial sum process $n^{-1/2} \sum_{k=1}^{[nt]}X_k$
  converges weakly to the Brownian motion? For $p\in[1,2)$, do the
  sequences $\{|X_n|^p\}$ have distributional long memory in the sense
  that $n^{-H} \sum_{k=1}^{[nt]}\{|X_k|^p-\mathbb{E}[|X_k|^p]\}$ converge to
  the fractional Brownian motion with Hurst index $H$ for a suitable
  $H>1/2$?
\item Statistical inference. The FIGARCH($p,d,q$) is a parametric model, so the
  issue of estimation of its parameter is naturally raised. Also, if $d$ is
  linked to some memory property of the process, semi-parametric estimation of
  $d$ would be of interest.
\end{enumerate}

\section*{Acknowledgement} The authors are grateful to Remigijus
Leipus for his comments on a preliminary version of this work.


\end{document}